\documentclass[12pt,fullpage]{article}

\usepackage{amsfonts}
\usepackage{xr}
\usepackage{graphicx}
\usepackage{amsmath}
\usepackage{epsfig}

\def\R{\mathbb{R}}

\def\epsilon{\varepsilon}

\def\tilde{\widetilde}
\def\div{\mbox{div}}

\def\trait (#1) (#2) (#3){\vrule width #1pt height #2pt depth #3pt}
\def\fin{\hfill\trait (0.1) (5) (0) \trait (5) (0.1) (0) \kern-5pt
\trait (5) (5) (-4.9) \trait (0.1) (5) (0)}

\newcommand{\SE}{\setcounter{equation}{0} \section}
\newcommand{\be}{\begin{equation}}
\newcommand{\ee}{\end{equation}}
\newcommand{\baa}{\begin{array}}
\newcommand{\eaa}{\end{array}}
\newcommand{\ba}{\begin{eqnarray}}
\newcommand{\ea}{\end{eqnarray}}

\newtheorem{theo}{\bf Theorem}[section]
\newtheorem{lem}[theo]{\bf Lemma}
\newtheorem{pro}[theo]{\bf Proposition}
\newtheorem{cor}[theo]{\bf Corollary}
\newtheorem{defi}[theo]{\bf Definition}
\newtheorem{rem}[theo]{\bf Remark}

\begin{document}
\date{}
\title{\bf{Divergence operator and Poincar\'e inequalities on arbitrary bounded domains}}
\author{Ricardo Duran$^{\hbox{\small{ a}}}$, Maria-Amelia Muschietti$^{\hbox{\small{ b}}}$, Emmanuel Russ$^{\hbox{\small{
c}}}$\\
and Philippe Tchamitchian$^{\hbox{\small{ c}}}$\\
\\
\footnotesize{$^{\hbox{a }}$ Universidad de Buenos Aires}\\
\footnotesize{Facultad de Ciencias Exactas y Naturales, Departmento de Matem\'atica}\\
\footnotesize{Ciudad Universitaria. Pabell\'on I, (1428) Buenos Aires, Argentina}\\
\footnotesize{$^{\hbox{b }}$ Universidad Nacional de La Plata}\\
\footnotesize{Facultad de Ciencias Exactas, Departamento de Matem\'atica}\\
\footnotesize{Casilla de Correo 172, 1900 La Plata, Provincia de Buenos Aires, Argentina}\\
\footnotesize{$^{\hbox{c }}$Universit\'e Paul C\'ezanne,CNRS,LATP (UMR 6632)}\\
\footnotesize{Facult\'e des Sciences et Techniques, LATP}\\
\footnotesize{Case cour A, Avenue Escadrille Normandie-Niemen, F-13397 Marseille
Cedex 20, France}}
\maketitle

\begin{abstract}
Let $\Omega$ be an arbitrary bounded domain of $\R^n$. We study the right invertibility of the divergence on $\Omega$ in weighted Lebesgue and Sobolev spaces on $\Omega$, and rely this invertibility to a geometric characterization of $\Omega$ and to weighted Poincar\'e inequalities on $\Omega$. We recover, in particular, well-known results on the right invertibility of the divergence in Sobolev spaces when $\Omega$ is Lipschitz or, more generally, when $\Omega$ is  a John domain, and focus on the case of $s$-John domains. 
\end{abstract}

\noindent{\small{{\bf AMS numbers 2000: }}} Primary, 35C15. Secondary, 35F05, 35F15, 46E35.


\tableofcontents


\SE{Introduction}\label{intro}
Let $n\geq 2$. Throughout this paper, $\Omega$ will denote a nonempty bounded domain ({\it i.e. } an open connected subset of $\R^n$), and, for all $p\in [1,+\infty]$, $L^p_0(\Omega)$ stands for the subspace of $L^p(\Omega)$ made of functions having zero integral on $\Omega$. \par
\noindent Let $1\leq p\leq +\infty$. Does there exist a constant $C>0$ such that, for all $f\in L^p_0(\Omega)$, one can find a vector field $u\in W^{1,p}_0(\Omega,\R^n)$ such that $\mbox{div }u=f$ in $\Omega$ and
$$
\left\Vert Du\right\Vert_{L^p(\Omega,\R^n)}\leq C\left\Vert f\right\Vert_{L^p(\Omega)} ?
$$
This problem arises in various contexts. In particular, it plays a key role (for $p=2$) in the analysis of the Stokes problem by finite elements methods (\cite{bs}, Chapter 10). It is also involved in the study of Hardy-Sobolev spaces on strongly Lipschitz domains of $\R^n$ (\cite{art}). \par
\noindent When $p=1$ or $p=+\infty$, the problem has no solution in general, even in $W^{1,p}(\Omega)$ (see \cite{dft} when $p=1$ and \cite{mcm} when $p=+\infty$). \par
\noindent When $1<p<+\infty$ and $\Omega$ is Lipschitz, the question receives a positive answer. A reformulation of this result is that the divergence operator, mapping $W^{1,p}_0(\Omega,\R^n)$ to $L^p_0(\Omega)$, has a continuous right inverse. \par
\noindent Several approaches can be found for this theorem in the literature. It can for instance be proved by functional analysis arguments (see \cite{art, necas}), or by a somewhat elementary construction (\cite{bb}). In dimension $2$ and when $\Omega$ is smooth (or is a convex polygon), it can be solved via the Neumann problem for the Laplace operator on $\Omega$ (\cite{ba, bs, gr, l}). This approach can also be adapted to the case when $\Omega$ is a non-convex polygon, relying on suitable trace theorems (\cite{asv}). \par

\medskip

\noindent A different and more explicit construction, valid in star-shaped domains with respect to a ball and for weighted norms, was developed in \cite{bog} and in \cite{dm}. The idea of this construction is to integrate test functions on segments included in $\Omega$, which yields a somewhat explicit formula for $u$. This approach was extended to John domains in \cite{adm}, replacing segments by appropriate rectifiable curves included in $\Omega$. \par

\medskip

\noindent However, there exist domains $\Omega$ such that $\mbox{div}:W^{1,p}_0(\Omega)\rightarrow L^p_0(\Omega)$ has no continuous right inverse. A counterexample, originally due to Friedrichs, can be found in \cite{adm} and is, actually, a planar domain with an external cusp of power type. In the present paper, we investigate the solvability of $\mbox{div }u=f$ on {\it arbitrary} bounded domains, with estimates in weighted Lebesgue or Sobolev spaces. \par
\noindent More precisely, we consider the following question: does there exist an integrable weight $\mbox{w}>0$ on $\Omega$ such that the divergence operator, acting from $L^{p}\left(\Omega,\frac 1{\mbox w}dx\right)$ to $L^p(\Omega,dx)$, has a continuous right inverse ? When $p=+\infty$, we give a necessary and sufficient condition on $\Omega$ for this property to hold. This condition says that the geodesic distance in $\Omega$ to a fixed point is integrable over $\Omega$. Under this assumption, the previous divergence problem is also solvable in $L^p$ norms for all $p\in (1,+\infty)$. \par
\noindent The proofs rely on two key tools: first, we extend the method developed in \cite{bog,dm,adm} to the case of an arbitrary bounded domain, using adapted rectifiable curves inside the domain. Then, we show that the solvability of the divergence problem in $L^p$ is equivalent to a weighted Poincar\'e inequality in $L^q$ with $\frac 1p+\frac 1q=1$.\par
\noindent When the divergence problem is solvable in $L^p$ spaces, we also solve it in weighted Sobolev spaces on $\Omega$. This extends to the case of arbitrary bounded domains the previously cited theorems when $\Omega$ is Lipschitz or, more generally, is a John domain. A key tool for this, interesting in itself, is a kind of atomic decomposition for functions in $L^p(\Omega)$ with zero integral. We also obtain interesting and new results for the invertibility of the divergence on Sobolev spaces in the class of $s$-John domains and, among these domains, in the sub-class of strongly h\"olderian domains. It should be noted that solvability results in weighted Sobolev spaces for some h\"olderian planar domains had already been obtained by the first author and F. Lopez Garcia in \cite{dg}. \par
\noindent In the case when $p=1$, as we said, it is impossible to solve $\mbox{div }u=f$ when $f\in L^1(\Omega)$ with $u\in W^{1,1}(\Omega)$ in general. However, if $f$ is supposed to be in a Hardy space on $\Omega$ instead of belonging to $L^1(\Omega)$, then, under suitable assumptions on $\Omega$, it is possible to solve in appropriate Hardy-Sobolev spaces, considered in \cite{art}. We will come back to this issue in a forthcoming paper.\par

\medskip

\noindent Here is the outline of the paper. In Section \ref{ineqlinfty}, we establish the equivalence between the solvability of the divergence in $L^{\infty}$ spaces and the integrabilty of the geodesic distance to a fixed point.  Section \ref{lpineq} is devoted to the solvability in $L^p$ spaces and the equivalence with $L^q$ Poincar\'e inequalities. In Section \ref{sobol}, we investigate the solvability of the divergence in weighted Sobolev spaces. Finally, in Section \ref{ex}, we focus on the case of $s$-John domains and, in particular, on the class of strongly h\"olderian domains.\par

\medskip

\noindent {\bf Acknowledgments: }
This work was supported by the bilateral CNRS/CONICET linkage (2007-2008): Equations aux d\'eriv\'ees partielles sur des domaines peu r\'eguliers / Partial Differential Equations on non-smooth domains.\par
\noindent The research was conducted at the University of Buenos Aires and the University of La Plata in Argentina, and at the Laboratoire d'Analyse, 
Topologie et Probabilit\'es at the Facult\'e des Sciences et Techniques de Saint-J\'er\^ome, Universit\'e Paul C\'ezanne, Marseille in France.\par
\noindent Some results of this work were presented by the first author in the fourth international symposium on nonlinear PDE's and free boundary problems, dedicated to Luis Caffarelli, in Mar del Plata, March 17th-20th 2009.

\SE{Invertibility of the divergence in $L^{\infty}$ spaces} \label{ineqlinfty}
Throughout this section, $\Omega$ denotes a bounded domain (connected open subset) of $\R^n$, without any assumption on the regularity of its boundary $\partial\Omega$. We investigate how to invert the divergence operator on $L^{\infty}$ spaces, asking the following question: does it exist a weight $\mbox{w}$, defined on $\Omega$, such that
\begin{itemize}
\item[$i)$]
$\mbox{w}$ is integrable over $\Omega$ ({\it{i.e.}} ${\mbox{w}} \in L^1(\Omega)$),
\item[$ii)$]
for each $f\in L^{\infty}(\Omega)$ with $\int_{\Omega} f(x)dx=0$, there exists a vector-valued function $u$ solution of
    \begin{equation} \label{divinftyw}
    \left\{
        \begin{array}{ll}
        \mbox{div }u=f &\mbox{ in } \Omega,\\
        u\cdot\nu=0 &\mbox{ on }\partial\Omega,
        \end{array}
    \right.
    \end{equation}
satisfying the estimate
    \begin{equation}\label{estimLinftyw}
    \left\Vert \frac{1}{\mbox{w}}u\right\Vert_{\infty} \leq C\left\Vert f\right\Vert_{\infty},
    \end{equation}
where $C>0$ only depends on $\Omega$ and the weight $\mbox{w}$?
\end{itemize}

\noindent We define $u$ being a solution of (\ref{divinftyw}) by demanding that
\begin{equation}\label{eqinE}
\int_{\Omega} u\cdot\nabla g = - \int_{\Omega} fg
\end{equation}
for all $g$ in the space $E$ of those functions in $L^1(\Omega)$ such that their distributional gradient is a function and
$$
\left\Vert u \right\Vert_{E} = \int_{\Omega} \left\vert g \right\vert + \int_{\Omega} \mbox{w} \left\vert \nabla g \right\vert < + \infty .
$$
The integrability condition $i)$ above ensures that ${\mathcal D}(\overline{\Omega})$, the set of restrictions to $\Omega$ of test functions in $\R^n$, is embedded in $E$. Thus, equation (\ref{eqinE}) defines the meaning of both the partial differential equation and the boundary condition in (\ref{divinftyw}).\par

\medskip
\noindent We say that the problem $(\div)_{\infty}$ is solvable in $\Omega$ when this question receives a positive answer.\par

\medskip
\noindent We introduce $\mbox{dist}_{\Omega}$, the euclidean geodesic distance in $\Omega$, defined by
$$
\mbox{dist}_{\Omega} (y,x) =  \inf l(\gamma),
$$
where $\gamma$ is any rectifiable path defined on $[0,1]$ such that $\gamma(0) = y$, $\gamma(1) = x$, and $\gamma(t) \in \Omega$ for all $t \in [0,1]$, $l(\gamma)$ being its length. We pick once and for all a reference point $x_0 \in \Omega$, and we simply denote by $d_{\Omega}(y)$ the number $\mbox{dist}_{\Omega}(y,x_0)$.\par

\medskip
\noindent Our main result is the following:
\begin{theo}\label{solvability}
The problem $(\div)_{\infty}$ is solvable in $\Omega$ if and only if $d_{\Omega} \in L^1(\Omega)$.
\end{theo}
Notice that $d_{\Omega}$ being integrable on $\Omega$ does not depend on the choice of the reference point $x_0$.\par

\medskip
\noindent To prove this theorem, assume first that $(\div)_{\infty}$ is solvable in $\Omega$. We will explain in the next section that it is equivalent to various Poincar\'e inequalities, implying in particular the following one:
$$
\int_{\Omega} \left\vert f(x) \right\vert dx \leq C \int_{\Omega} {\mbox{w}}(x) \left\vert \nabla f(x) \right\vert dx
$$
for all $f \in {\mbox{W}}^{1,1}_{loc}(\Omega)$ which vanishes on a non negligible subset $K$ of $\Omega$, the constant $C$ depending on $\Omega$ and on the ratio $\frac{\left\vert \Omega \right\vert}{\left\vert K \right\vert}$.\par
\noindent We apply this inequality to
$$
f(x) = \left( d_{\Omega}(x) - r_0 \right)_{+},
$$
where $r_0$ is the radius of some fixed ball centered at $x_0$ and included in $\Omega$. Since
$$
\left\vert d_{\Omega}(y) - d_{\Omega}(x) \right\vert \leq \left\vert y-x \right\vert
$$
for all close enough $x, y \in \Omega$, we have $\left\vert \nabla f(x) \right\vert \leq 1$ almost everywhere, so that
$$
\int_{\Omega} {\mbox{w}}(x) \left\vert \nabla f(x) \right\vert dx \leq \int_{\Omega} {\mbox{w}}(x) dx < +\infty .
$$
Thus $f \in L^{1}(\Omega)$, which implies that $d_\Omega$ is integrable on $\Omega$, too.\par

\medskip
\noindent Reciprocally we assume that
\begin{equation}\label{dOmega}
\int_{\Omega} d_{\Omega}(y) dy < +\infty ,
\end{equation}
and explain how to solve $(\div)_{\infty}$ under this hypothesis. Our method will be constructive.\par

\medskip
\noindent The starting point is a family of curves relying each point in $\Omega$ to $x_0$. For the sake of simplicity, the image of any path $\gamma$ will be called $\gamma$, too.\par
\noindent We use a Whitney decomposition of $\Omega$, whose properties we recall. It is, up to negligible sets, a partition of $\Omega$ in closed cubes:
\[
\Omega=\bigcup_{j\geq 0} Q_j,
\]
where the interiors of the $Q_j$'s are pairwise disjoint, and for all $j\geq 0$, $2Q_j\subset \Omega$ while $4Q_j\cap \Omega^c\neq \emptyset$. We denote by $x_j$ and $l_j$ the center and the side length of $Q_j$ (we may and do assume that $x_0$ is the center of $Q_0$). A key observation of constant use is that, if $x \in \Omega$, then
$$
\frac12 l_j \leq d(x) \leq \frac{5 \sqrt n}{2} l_j
$$
whenever $x \in Q_j$. We denote by $d(x)$ the distance from $x$ to the boundary of $\Omega$.\par
\noindent We select a path $\gamma_j : [0,1] \mapsto \Omega$ such that $\gamma_{j}(0) = x_j$, $\gamma_{j}(1) = x_0$, $l(\gamma_j) \leq 2 d_{\Omega}(x_j)$. When $\gamma_j$ intersects some Whitney cube $Q_k$, we replace $\gamma_j \cap Q_k$ by the segment joining its two endpoints. In particular, we call $[x_j, \tilde{x}_j]$ the segment $\gamma_j \cap Q_j$. The modified path is still denoted by $\gamma_j$. We certainly have not increased its length, and moreover we have ensured that
$$
l(\gamma_j \cap B(x, r)) \leq C r
$$
for every $x \in \Omega$ and $r \leq \frac12 d(x)$, where $C$ only depends on the dimension.  \par
\noindent Finally, if $y \in Q_j$ for some $j$ (when there are several such $j$'s, just arbitrarily select one of them), we link $y$ to $\tilde{x}_j$ by a segment, and then $\tilde{x}_j$ to $x_0$ by $\gamma_j$. The resulting path linking $y$ to $x_0$ is called $\gamma(y)$, and its current point $\gamma(t,y)$, where $t\in [0,1]$. By construction, the following properties hold:
\begin{itemize}
\item[$(\gamma.a)$]
for all $y\in \Omega$ and $t \in [0,1]$, $\gamma(t,y)\in \Omega$,
with
$\gamma(0,y)=y$, $\gamma(1,y)=x_0$,
\item[$(\gamma.b)$]
$(t,y)\mapsto \gamma(t,y)$ is measurable and $\gamma(y)$ is rectifiable,
\item[$(\gamma.c)$]
there exists a constant $C > 0$, only depending on the dimension, such that for all $x, y \in \Omega$
    \begin{equation} \label{Ahlfors}
    \forall r \leq \frac 12 d(x), \ l(\gamma(y) \cap B(x, r)) \leq C r
    \end{equation}
and
	\begin{equation} \label{length}
	l(\gamma(y)) \leq C d_{\Omega}(y),
	\end{equation}
\item[$(\gamma.d)$]
for all $\varepsilon >0$ small enough, there exists $\delta >0$ such that
$$
\forall y \in \Omega_{\varepsilon}, \ \gamma(y) \subset \Omega_{\delta},
$$
where by definition $\Omega_{s} = \{ z\in \Omega ; d(z) > s \}$.
\end{itemize} \par

\medskip
\noindent We now define a weight $\omega$ on $\Omega$ (implicitly depending on the family of paths $\gamma$) by
\begin{equation} \label{defomega}
\omega(x)=\left\vert \left\{y\in \Omega;\ \mbox{dist}_{\Omega}(\gamma(y),x)\leq \frac 12 d(x)\right\}\right\vert,
\end{equation}
where $\mbox{dist}_{\Omega}(\gamma(y),x)$ is the distance in $\Omega$ from the path $\gamma(y)$ to the point $x$. We first prove that:
\begin{lem} \label{omega}
 \begin{equation}\label{intomega}
    \int_\Omega \omega(x) d(x)^{-n+1} dx < + \infty.
    \end{equation}
\end{lem}

\noindent {\bf Proof:} Let $y$ be such that there exists $t_0$ satisfying $\mbox{dist}_{\Omega}(\gamma(t_0,y),x) \leq \frac1{2}d(x)$, which here is the same as $\left\vert \gamma(t_0,y) - x \right\vert \leq \frac1{2}d(x)$. Then we have $\left\vert \gamma(t,y) - x \right\vert < \frac{3}{4}d(x)$ whenever $\left\vert \gamma(t,y) - \gamma(t_0,y) \right\vert < \frac1{4} d(x)$. This implies that
$$
\begin{array}{ll}
d(x) & \leq 4 l(\gamma(y) \cap B(\gamma(t_0,y),\frac1{4}d(x)) \\
        & \leq 4 \int_{0}^{1} \left\vert \dot\gamma(t,y) \right\vert {\bf{1}}_{\left\vert \gamma(t,y) - x \right\vert < \frac{3}{4}d(x)} dt ,
\end{array}
$$
and therefore that
$$
\int_{\Omega} \omega(x) d(x)^{-n+1} dx \leq C \int_{\Omega} \int_{\Omega} \int_0^1 \left\vert \dot{\gamma}(t,y) \right\vert {\bf{1}}_{\left\vert x-\gamma(t,y) \right\vert < \frac{3}{4}d(x)} d(x)^{-n} dt dy dx.
$$
In the integral above, $d(x)$ is comparable with $d(\gamma(t,y))$ uniformly in $x$, $t$ and $y$. Integrating first with respect to $x$, this observation leads to
$$
\begin{array}{ll}
\int_{\Omega} \omega(x) d(x)^{-n+1} dx & \leq C \int_{\Omega} \int_0^1 \left\vert \dot{\gamma}(t,y) \right\vert dt dy \\
                                                                       & = C \int_{\Omega} l(\gamma(y)) dy ,
\end{array}
$$
and the proof is ended thanks to the hypothesis (\ref{dOmega}) and the estimate (\ref{length}).

\medskip

\noindent The key result is then the following, that we state in full generality for further use:
\begin{lem} \label{solutioninfty}
Assume that $\Omega$ fulfills(\ref{dOmega}). Let $\gamma = \left\{ \gamma(y), y \in \Omega \right\}$ be any family of paths satisfying properties $(\gamma.a)$ to $(\gamma.d)$, and $\omega$ defined by (\ref{defomega}). For each $f\in L^{\infty}(\Omega)$ with $\int_{\Omega} f(x)dx=0$, there exists a vector-valued function $u$ solution of
    \begin{equation} \label{divinfty}
    \left\{
        \begin{array}{ll}
        \mbox{div }u=f &\mbox{ in }\Omega,\\
        u\cdot\nu=0 &\mbox{ on }\partial\Omega,
        \end{array}
    \right.
    \end{equation}
satisfying the estimate
    \begin{equation}\label{estimLinfty}
    \left\Vert \frac{d^{n-1}}{\omega}u\right\Vert_{\infty} \leq C\left\Vert f\right\Vert_{\infty},
    \end{equation}
where $C>0$ only depends on $\Omega$ and the choice of the family $\gamma$.
\end{lem}

{\bf Proof: } it relies on a representation formula, which goes back (in its earlier version) to \cite{bog}. A first generalization, applicable to John domains, was designed in \cite{adm}. Here, we still generalize it to the case of an arbitrary domain.\par

\medskip By dilating and translating $\Omega$, we may and do assume that $x_0=0$ and $d(0)=15$. We choose a function $\chi\in {\mathcal D}(\Omega)$ supported in $\overline{B(0,1)}$ and such that $\int_{\Omega} \chi(x)dx=1$. For each $y\in \Omega$, let $\tau(y)$ be the smallest $t>0$ such that $\gamma(t,y) \in \partial B(y, \frac12 d(y))$ (if there is no such $t$, just choose $\tau(y) = 1$). We define a function $t\mapsto \rho(t,y)$, $t\in [0,1]$, by:
$$
\begin{array}{ll}
\rho(t,y)=\alpha \left\vert y - \gamma(t,y) \right\vert &\mbox{ if }t\leq \tau(y),\\
\rho(t,y)=\frac 1{15} d(\gamma(t,y)) &\mbox{ if }t>\tau(y),
\end{array}
$$
where $\alpha$ is so chosen that $\rho(\cdot, y)$ is a continuous function.
This means that
$$
\alpha = \frac{2}{15} \frac{d(\gamma(\tau(y),y))}{d(y)},
$$
and the reader may check that $\alpha \leq \frac 15$.
By construction, we always have
\begin{equation}\label{a}
\rho(t,y)\leq \frac 1{5}d(\gamma(t,y))
\end{equation}
so that
$$
\gamma(t,y)+\rho(t,y)z \in \Omega
$$
for every $t \in [0,1]$ and $z \in B(0,1)$. This comes from the fact that, if $t \leq \tau(y)$, then $\left\vert y - \gamma(t,y) \right\vert \leq \frac12 d(y)$, which implies $\rho(t,y) \leq \frac{\alpha}{2} d(y)$ and $d(y) \leq 2 d(\gamma(t,y))$.

\medskip
\noindent Let us start with $\varphi \in {\mathcal D}(\R^n)$, $y \in \Omega$ and $z \in B(0,1)$. We have
$$
\varphi (y) - \varphi (z) = - \int_0^1 (\dot{\gamma}(t,y) + \dot{\rho}(t,y)z) \cdot \nabla\varphi(\gamma(t,y) + \rho(t,y)z) dt.
$$
Multiplying by $\chi (z)$ and integrating, we get
$$
\varphi (y) - \int_{\Omega}\varphi\chi = - \int_{B(0,1)}\int_0^1 (\dot{\gamma}(t,y) + \dot{\rho}(t,y)z) \cdot \nabla\varphi(\gamma(t,y) + \rho(t,y)z) \chi(z) dtdz.
$$
Since $\int_{\Omega}f = 0$, this implies
$$
\int_{\Omega}f\varphi = - \int_{\Omega}\int_{B(0,1)}\int_0^1 f(y) (\dot{\gamma}(t,y) + \dot{\rho}(t,y)z) \cdot \nabla\varphi(\gamma(t,y) + \rho(t,y)z) \chi(z) dtdzdy.
$$
Changing $z$ into $x = \gamma(t,y) + \rho(t,y)z$, this formula becomes
\begin{equation}\label{fphi}
\int_{\Omega}f\varphi = - \int_{\Omega}\int_{\Omega}\int_0^1 f(y) \left(\dot{\gamma}(t,y) + \dot{\rho}(t,y)\frac{x-\gamma(t,y)}{\rho(t,y)}\right) \cdot \nabla\varphi(x) \chi\left(\frac{x-\gamma(t,y)}{\rho(t,y)}\right) \frac1{\rho(t,y)^n}dtdxdy.
\end{equation}
If $x, y \in \Omega$, we define the vector-valued kernel $G$ by
$$
G(x,y) = \int_0^1 \left[\dot{\gamma}(t,y) + \dot{\rho}(t,y)\frac{x-\gamma(t,y)}{\rho(t,y)}\right] \chi\left(\frac{x-\gamma(t,y)}{\rho(t,y)}\right) \frac{dt}{\rho(t,y)^n} .
$$
Thanks to the support condition on $\chi$, we must have
$$
\left\vert x - \gamma(t,y) \right\vert < \rho(t,y)
$$
for the integrated term to be non zero. Hence $G(x,y)$ is well defined as soon as $x \neq y$. We moreover have
\begin{lem}\label{lemG}
There exists a constant $C>0$, depending on $\Omega$ and on the paths $\gamma$ only, such that
    \begin{equation}\label{G}
    \forall x \in \Omega, \   \int_\Omega \left\vert G(x,y) \right\vert dy \leq C \omega(x) d(x)^{-n+1}.
    \end{equation}
\end{lem}
Let us admit this statement for the moment and finish the proof of Lemma \ref{solutioninfty}. We define $u$ by
$$
u(x) = \int_\Omega G(x,y) f(y) dy.
$$
The estimate (\ref{G}) shows that this is meaningful and that
\begin{equation}\label{u}
\left\vert u(x) \right\vert \leq C \left\Vert f \right\Vert_{\infty} \omega(x) d(x)^{-n+1}.
\end{equation}
Thus (\ref{estimLinfty}) is satisfied. Then, it follows from (\ref{fphi}) and Fubini theorem, which we may apply thanks to Lemma \ref{omega} and (\ref{G}), that
\begin{equation}\label{eqinDOmega}
\int_{\Omega} u\cdot\nabla \varphi = - \int_{\Omega} f\varphi ;
\end{equation}
this is (\ref{eqinE}) for $\varphi \in {\mathcal D}(\overline{\Omega})$.\par

\noindent Let us now take $g \in E$. If $\varepsilon >0$ is small enough, we define
$$
f_{\varepsilon} = f {\bf{1}}_{\Omega_{\varepsilon}}
$$
and
$$
u_{\varepsilon}(x) = \int_\Omega G(x,y) f_{\varepsilon}(y) dy.
$$
We have $u_{\varepsilon}(x) \rightarrow u(x)$ for all $x \in \Omega$, and
\begin{equation}\label{uepsilon}
\left\vert u_{\varepsilon}(x) \right\vert \leq C \left\Vert f \right\Vert_{\infty} \omega(x) d(x)^{-n+1},
\end{equation}
by Lemma \ref{lemG}. Because of property $(\gamma.d)$ there exists $\delta >0$ such that ${\mbox{Supp }}u_{\varepsilon} \subset \Omega_{\delta}$. Moreover, there exists $C(\delta)$ satisfying
$$
\forall x \in \Omega_{\delta} \ \frac{1}{C(\delta)} \leq \omega (x) d(x)^{-n+1} \leq C(\delta).
$$
Hence, we may use a standard approximation argument to deduce from (\ref{eqinDOmega}), applied to $f_\varepsilon$ and $u_\varepsilon$, that
$$
\int_{\Omega} u_{\varepsilon}\cdot\nabla g = - \int_{\Omega} f_{\varepsilon} g.
$$
We conclude by using (\ref{uepsilon}) and the dominated convergence theorem to obtain:
$$
\int_{\Omega} u\cdot\nabla g = - \int_{\Omega} f g
$$
for all $g \in E$.\par

\medskip
\noindent We turn to the proof of Lemma \ref{lemG}. Since $\rho$ is defined according to two different cases, we decompose
$$
G(x,y) = G_1 (x,y) + G_2 (x,y)
$$
with
$$
G_1 (x,y) = \int_0^{\tau(y)} \left[\dot{\gamma}(t,y) +\alpha \frac{\dot{\gamma}(t,y) \cdot (y-\gamma(t,y))}{ \left\vert y - \gamma(t,y) \right\vert}\right] \chi\left(\frac{x-\gamma(t,y)}{\alpha \left\vert y - \gamma(t,y) \right\vert}\right) \frac{1}{\alpha^n  \left\vert y - \gamma(t,y) \right\vert^n} dt,
$$
and
\begin{equation}\label{G2}
G_2 (x,y) = \int_{\tau(y)}^1 \left[\dot{\gamma}(t,y) + [\dot{\gamma}(t,y) \cdot \nabla d(\gamma(t,y))] \frac{x-\gamma(t,y)}{d(\gamma(t,y))} \right] \chi\left(15\frac{x-\gamma(t,y)}{d(\gamma(t,y))}\right) \frac{15^n}{d(\gamma(t,y))^n} dt.
\end{equation}
We first estimate $\int_\Omega \left\vert G_1 (x,y) \right\vert dy$. When $\chi(\frac{x-\gamma(t,y)}{\alpha \left\vert y - \gamma(t,y) \right\vert})$ is non zero, we have
\begin{equation} \label{alpha}
\begin{array}{ll}
\left\vert x-y \right\vert & \leq (1+\alpha) \left\vert y-\gamma(t,y) \right\vert \\
                                       & \leq \frac{1+\alpha}{2} d(y) \\
                                       & \leq \frac{1+\alpha}{1-\alpha} d(x)
\end{array}
\end{equation}
and
$$
\begin{array}{ll}
\left\vert x - \gamma(t,y) \right\vert  &  \leq \alpha \left\vert y - \gamma(t,y) \right\vert \\
                                                             &  \leq \frac{\alpha}{1-\alpha} \left\vert x - y \right\vert.
\end{array}
$$
Therefore we obtain
$$
\begin{array}{ll}
\left\vert G_1 (x,y) \right\vert & \leq C \int_0^{\tau(y)} \left\vert \dot{\gamma}(t,y) \right\vert {\bf{1}}_{\left\vert x - \gamma(t,y) \right\vert \leq \frac{\alpha}{1-\alpha} \left\vert x - y \right\vert}  dt  \left\vert x-y\right\vert^{-n} \\
                                                  & \leq C  \left\vert x-y\right\vert^{-n} l(\gamma(y) \cap B(x, \frac{\alpha}{1-\alpha} \left\vert x - y \right\vert)).
\end{array}
$$
Since $\alpha \leq \frac15$, we deduce from (\ref{alpha}) that $\frac{\alpha}{1-\alpha} \left\vert x - y \right\vert \leq \frac38 d(x)$. We may therefore apply property $(\gamma.c)$, which gives
$$
\left\vert G_1 (x,y) \right\vert \leq C \left\vert x - y \right\vert^{-n+1}.
$$
Applying (\ref{alpha}) again, we see that
$$
\int_\Omega \left\vert G_1 (x,y) \right\vert dy \leq C d(x).
$$
By definition of $\omega$ there is a constant $c$ depending on the dimension such that $\omega(x) \geq c d(x)^n$, since $B(x, \frac1{2}d(x)) \subset \{y ; \mbox{dist}_{\Omega}(\gamma(y),x) \leq \frac1{2}d(x)\}$. Thus the analog of inequality (\ref{G}) is satisfied by $G_1$.\\
We now estimate $\int_\Omega \left\vert G_2 (x,y) \right\vert dy$. We have from (\ref{G2}):
$$
\left\vert G_2 (x,y) \right\vert \leq C \int_{\tau(y)}^1 \left\vert \dot{\gamma}(t,y) \right\vert {\bf{1}}_{\left\vert x-\gamma(t,y) \right\vert < \frac1{15}d(\gamma(t,y))} \frac{dt}{d(\gamma(t,y))^n}.
$$
On the domain of integration, the ratio $\frac{d(x)}{d(\gamma(t,y))}$ is bounded from above and below by $16/15$ and $14/15$. This implies the estimate
$$
\left\vert G_2 (x,y) \right\vert \leq C d(x)^{-n} l(\gamma(y) \cap B(x,\frac1{14}d(x))),
$$
so that, by property $(\gamma.c)$, we have
$$
\left\vert G_2 (x,y) \right\vert \leq C d(x)^{-n+1}
$$
uniformly in $x$, $y$, and $G_2(x,\cdot)$ is supported in the set $\{y ; d(\gamma(y),x) < \frac1{14}d(x)\}$. We therefore obtain the analog of inequality (\ref{G}) for $G_2$ as well, concluding the proof of Lemma \ref{lemG}.\\
\SE{Invertibility of the divergence in $L^p$ spaces} \label{lpineq}
The present section is devoted to the proof of the following $L^p$-type result:
\begin{theo} \label{divp}
Let $p\in (1,+\infty]$ and $\Omega\subset \R^n$ be an arbitrary bounded domain. Assume that the function $d_{\Omega}$ is integrable on $\Omega$. 
Then, if $f\in L^p(\Omega)$ with $\int_{\Omega} f(x)dx=0$, there exists a vector-valued function $u$ solution of
\begin{equation} \label{divlp}
\left\{
\begin{array}{ll}
\mbox{div }u=f &\mbox{ in }{\mathcal D}^{\prime}(\Omega),\\
u\cdot\nu=0 &\mbox{ on }\partial\Omega,
\end{array}
\right.
\end{equation}
satisfying the estimate
\begin{equation} \label{estimlp}
\left\Vert \frac{d^{n-1}}{\omega}u\right\Vert_{p} \leq C\left\Vert f\right\Vert_{p},
\end{equation}
where $C>0$ only depends on $\Omega$ and the choice of the paths $\gamma$. 
\end{theo}
Let us clarify the meaning of (\ref{divlp}). Define $q\in [1,+\infty)$ by $\frac 1p+\frac 1q=1$ and let
$$
E^q(\Omega):=\left\{g\in L^q\left(\Omega\right);\ \int_{\Omega} \left\vert \nabla g(x)\right\vert^q \left(\frac{\omega(x)}{d^{n-1}(x)}\right)^q<+\infty\right\},
$$
equipped with the norm
$$
\left\Vert g\right\Vert_{E^q(\Omega)}:=\left\Vert g\right\Vert_{L^q(\Omega)}+ \left\Vert \frac{\omega(x)}{d^{n-1}(x)} \left\vert \nabla g\right\vert\right\Vert_{L^q(\Omega)}.
$$
By (\ref{divlp}), we mean that 
$$
\int_{\Omega} u\cdot \nabla \varphi=\int_{\Omega} f\varphi
$$
for all $\varphi\in E^q(\Omega)$.\par

\medskip

\noindent The proof of Theorem \ref{divp} goes through a family of Poincar\'e inequalities, which we present now. 
If $1\leq q<+\infty$, say that $\Omega$ supports a weighted $L^q$ Poincar\'e inequality if and only if, for all function $g\in E^q(\Omega)$,
\begin{equation} \label{Pq} \tag{$P_q$}
\left\Vert g-g_{\Omega}\right\Vert_{L^q(\Omega)}\leq C\left\Vert \frac{\omega}{d^{n-1}}\left\vert \nabla g\right\vert\right\Vert_{L^q(\Omega)},
\end{equation}
where
$$
g_{\Omega}:=\frac 1{\left\vert \Omega\right\vert} \int_{\Omega} g(x)dx.
$$
The strategy of the proof of Theorem \ref{divp} is as follows: we first establish the general fact that, for all $p\in (1,+\infty]$, the solvability of (\ref{divlp}) for all $f\in L^p(\Omega)$ having zero integral with the estimate (\ref{estimlp}) is equivalent to the validity of (\ref{Pq}), with $\frac 1p+\frac 1q=1$. We then prove that $(P_1)$ implies $(P_q)$ for any $q\in [1,+\infty)$. It is then a consequence of Lemma \ref{solutioninfty} that, since $d_{\Omega}\in L^1(\Omega),$ $(P_1)$, and therefore $(P_q)$, hold for all $q\in [1,+\infty)$ and Theorem \ref{divlp} therefore follows again from the equivalence with $(P_q)$.\par
\noindent In the sequel, we will say that $\left(\mbox{div}_p\right)$ holds if and only if, for all $f\in L^p(\Omega)$ with zero integral,  there exists a vector-valued function $u$ solution of (\ref{divlp}) such that the estimate (\ref{estimlp}) holds with a constant $C>0$ only depending on $\Omega$ and the choice of the paths $\gamma$. 
\subsection{Solvability of the divergence and Poincar\'e inequalities}
Let us prove the following equivalence:
\begin{pro} \label{equivdivpoinc}
Let $\Omega\subset \R^n$ be a  bounded domain. Let $1<p\leq +\infty$ and $1\leq q<+\infty$ be such that $\frac 1p+\frac 1q=1$. Then $(P_q)$ holds if and only if $\left(\mbox{div}_p\right)$ holds.
\end{pro}
{\bf Proof: } Assume first that $\left(\mbox{div}_p\right)$ holds. Let $f\in E^q(\Omega)$ and $g\in L^p(\Omega)$ with $\left\Vert g\right\Vert_p\leq 1$. Pick up a solution $u$ of 
$$
\left\{
\begin{array}{ll}
\displaystyle \mbox{div }u=g-g_{\Omega} & \mbox{ in }\Omega,\\
\displaystyle u\cdot \nu=0 &\mbox{ on }\partial\Omega
\end{array}
\right.
$$
satisfying
$$
\left\Vert \frac{d^{n-1}}{\omega}u\right\Vert_p\leq C\left\Vert g\right\Vert_p.
$$
Then, by definition of $u$,
\[
\begin{array}{lll}
\displaystyle \left\vert \int_{\Omega} \left(f(x)-f_{\Omega}\right)g(x)dx\right\vert & =& \displaystyle \left\vert \int_{\Omega} \left(f(x)-f_{\Omega}\right)\left(g(x)-g_{\Omega}\right)dx\right\vert\\
& = & \displaystyle \left\vert \int_{\Omega} \nabla f(x)\cdot u(x)dx\right\vert \\
& \leq & \displaystyle \left\Vert \frac{\omega}{d^{n-1}} \left\vert \nabla f\right\vert \right\Vert_{q} \left\Vert \frac{d^{n-1}}{\omega} u\right\Vert_{p}\\
& \leq & \displaystyle C \left\Vert \frac{\omega}{d^{n-1}} \left\vert \nabla f\right\vert \right\Vert_{q} \left\Vert g\right\Vert_p\\
& \leq & \displaystyle C \left\Vert \frac{\omega}{d^{n-1}} \left\vert \nabla f\right\vert \right\Vert_{q},
\end{array}
\]
which proves that 
$$
\left\Vert f-f_{\Omega}\right\Vert_q\leq C\left\Vert \frac{\omega}{d^{n-1}} \left\vert \nabla f\right\vert \right\Vert_{q},
$$
and this means that $(P_q)$ holds.\par
\medskip

\noindent Assume conversely that $(P_q)$ holds and let $f\in L^p(\Omega)$ with zero integral. Define
$$
{\mathcal G}_q:=\left\{v\in L^q\left(\Omega,\R^n,\frac{\omega}{d^{n-1}}\right);\ v=\nabla g\mbox{ for some }g\in E^q(\Omega)\right\},
$$
equipped with the norm of $L^q\left(\Omega,\R^n,\frac{\omega}{d^{n-1}}\right)$.\par
\noindent Define a linear form on ${\mathcal G}_q$ by
$$
Lv=\int_{\Omega} f(x)g(x)dx\mbox{ if }v=\nabla g.
$$
Observe that $L$ is well-defined since $f$ has zero integral on $\Omega$ and that
$$
\begin{array}{lll}
\displaystyle \left\vert L(v)\right\vert &= & \displaystyle \left\vert \int_{\Omega} f(x)\left(g(x)-g_{\Omega}\right)dx\right\vert\\
& \leq & \displaystyle C\left\Vert f\right\Vert_p \left\Vert \frac{\omega}{d^{n-1}}\left\vert \nabla g\right\vert\right\Vert_q
\end{array}
$$
since $(P_q)$ holds; this shows that $L$ is bounded on ${\mathcal G}_q$. By the Hahn-Banach theorem, $L$ may be extended to a bounded linear form on $L^q\left(\Omega,\R^n,\frac{\omega}{d^{n-1}}\right)$ with functional norm bounded by $C\left\Vert f\right\Vert_p$. Therefore, since $1\leq q<+\infty$, there exists a vector field $u\in L^p\left(\Omega,\R^n,\frac{\omega}{d^{n-1}}\right)$ such that, for all $g\in E^q(\Omega)$,
\begin{equation} \label{dualug}
\int_{\Omega} u(x)\cdot \nabla g(x)dx=\int_{\Omega} f(x)g(x)dx,
\end{equation}
with
\begin{equation} \label{estimup}
\left\Vert \frac{d^{n-1}}{\omega}u\right\Vert_{p}\leq C\left\Vert f\right\Vert_p.
\end{equation}
Note that, for $p=1$, this uses the fact that the measure $\frac{\omega}{d^{n-1}}dx$ is finite (Lemma \ref{omega}). The identity (\ref{dualug}) means that $u$ solves (\ref{divlp}), and (\ref{estimup}) is exactly (\ref{estimlp}). This ends the proof of Proposition \ref{equivdivpoinc}.\hfill\fin\par

\subsection{From an $L^1$ to an $L^q$ Poincar\'e inequality}
We now prove the following fact:
\begin{pro} \label{poinc1q}
Assume that $(P_1)$ holds. Then, for all $q\in [1,+\infty[$, $(P_q)$ holds.
\end{pro}
The proof mimics a usual procedure for Poincar\'e inequalities (see \cite{hk}) and we give it for the sake of completeness. It relies on a very useful characterization of Poincar\'e inequalities:
\begin{lem} \label{poincarecarac}
Let $1\leq q<+\infty$. Then, $(P_q)$ holds if and only if there exists $C>0$ with the following property: for all measurable subset $E\subset \Omega$ with $\left\vert E\right\vert\geq \frac 12\left\vert \Omega\right\vert$ and for all $g\in E^q(\Omega)$ vanishing on $E$, 
$$
\left\Vert g\right\Vert_q\leq C\left\Vert \frac{\omega}{d^{n-1}}\left\vert \nabla g\right\vert\right\Vert_q.
$$
\end{lem}


\noindent {\bf Proof of Lemma: } assume first that $(P_q)$ holds and let $E$ and $g$ be as in the statement of the lemma. Notice that
$$
\begin{array}{lll}
\displaystyle \left\vert E\right\vert \left\vert g_{\Omega}\right\vert^q &= & \displaystyle \int_E \left\vert g(x)-g_{\Omega}\right\vert^qdx\\
& \leq & \displaystyle  \int_{\Omega} \left\vert g(x)-g_{\Omega}\right\vert^qdx\\
& \leq & \displaystyle C\int_{\Omega} \left\vert \frac{\omega(x)}{d^{n-1}(x)} \nabla g(x)\right\vert^q dx,
\end{array}
$$
which shows that
$$
\left\vert g_{\Omega}\right\vert\leq C^{\prime}\left(\frac 1{\left\vert \Omega\right\vert} \int_{\Omega} \left\vert \frac{\omega(x)}{d^{n-1}(x)} \nabla g(x)\right\vert^q dx\right)^{\frac 1q}.
$$
Since $\left\Vert g\right\Vert_q\leq \left\Vert g-g_{\Omega}\right\Vert_q+ \left\vert \Omega\right\vert^{\frac 1q} \left\vert g_{\Omega}\right\vert$, one gets the desired conclusion. \par

\medskip

\noindent For the converse, take $g\in L^q(\Omega)$. Observe first that there exists $\lambda\in \R$ such that, if
$$
E_{\lambda}:=\left\{x\in \Omega;\ g(x)\geq \lambda\right\}\mbox{ and } F_{\lambda}:= \left\{x\in \Omega;\ g(x)\leq \lambda\right\},
$$
then
$$
\left\vert E_{\lambda}\right\vert \geq \frac 12 \left\vert \Omega\right\vert\mbox{ and } \left\vert F_{\lambda}\right\vert\geq \frac 12\left\vert \Omega\right\vert.
$$
Indeed, for all $\lambda\in \R$, set
$$
\mu(\lambda):=\left\{x\in \Omega;\ g(x)\leq \lambda\right\}.
$$
The function $\mu$ is non-decreasing, right-continuous and satisfies 
$$
\lim_{\lambda\rightarrow -\infty} \mu(\lambda)=0\mbox{ and } \lim_{\lambda\rightarrow +\infty} \mu(\lambda)=\left\vert \Omega\right\vert.
$$
Therefore, $\lambda:=\inf\left\{t\in \R;\ \mu(t)\geq \frac 12\left\vert \Omega\right\vert\right\}$ satisfies the required properties. \par
\noindent Observe that $(g-\lambda)_+\in E^q(\Omega)$ (this can be proved by approximation arguments analogous to those used in the theory of usual Sobolev spaces). Since $(g-\lambda)_+$ vanishes on $F_{\lambda}$, the assumption yields
$$
\left\Vert (g-\lambda)_+\right\Vert_q\leq C\left\Vert \frac{\omega}{d^{n-1}}\left\vert \nabla g\right\vert\right\Vert_q.
$$
Similarly, since $(g-\lambda)_-$ vanishes on $E_{\lambda}$,
$$
\left\Vert (g-\lambda)_-\right\Vert_q\leq C\left\Vert \frac{\omega}{d^{n-1}}\left\vert \nabla g\right\vert\right\Vert_q,
$$
and we therefore obtain
$$
\left\Vert g-\lambda\right\Vert_q\leq C\left\Vert \frac{\omega}{d^{n-1}}\left\vert \nabla g\right\vert\right\Vert_q,
$$
from which $(P_q)$ readily follows. \hfill\fin\par

\begin{rem}
Lemma \ref{poincarecarac} extends straightforwardly  to the case when $E$ is any non-negligible measurable subset of $\Omega$, the constant $C$ depending on the ratio $\frac{\left\vert \Omega\right\vert}{\left\vert E\right\vert}$.
\end{rem}
\medskip

\noindent{\bf Proof of Proposition \ref{poinc1q}: } assume that $(P_1)$ holds, let $q\in [1,+\infty[$ and $g\in E^q(\Omega)$ vanishing on a subset $E\subset \Omega$ satisfying $\left\vert E\right\vert\geq \frac 12\left\vert \Omega\right\vert$. Again, it can be proved as for usual Sobolev spaces that $\left\vert g\right\vert^q\in E_1(\Omega)$. Applying then $(P_1)$ to $\left\vert g\right\vert^q$, which also vanishes on $E$, and using Lemma \ref{poincarecarac}, we obtain
$$
\begin{array}{lll}
\displaystyle \int_{\Omega} \left\vert g(x)\right\vert^qdx & \leq & \displaystyle C\int_{\Omega} \frac{\omega(x)}{d^{n-1}(x)} \left\vert \nabla \left(\left\vert g\right\vert^q\right)\right\vert (x)dx\\
& = & \displaystyle Cq\int_{\Omega} \frac{\omega(x)}{d^{n-1}(x)} \left\vert g(x)\right\vert^{q-1}\left\vert \nabla g(x)\right\vert dx\\
& \leq & \displaystyle Cq \left(\int_{\Omega}   \left\vert g(x)\right\vert^qdx\right)^{1-\frac 1q} \left(\int_{\Omega} \left(\frac{\omega(x)}{d^{n-1}(x)}\right)^q \left\vert \nabla g(x)\right\vert^qdx\right)^{\frac 1q},
\end{array}
$$
which yields exactly $(P_q)$ by Lemma \ref{poincarecarac} again. \hfill\fin\par

\SE{Solvability of the divergence in weighted Sobolev spaces} \label{sobol}
We now solve the divergence equation in weighted Sobolev spaces. For $p\in (1,+\infty)$, we define
$$
W^{1,p}\left(\Omega,\frac{d^n}{\omega}\right):=\left\{g\in L^p\left(\Omega,\frac{d^n}{\omega}\right);\ \partial_ig\in L^p\left(\Omega,\frac{d^n}{\omega}\right)\mbox{ for all }1\leq i\leq n\right\}
$$
and $W^{1,p}_0\left(\Omega,\frac{d^n}{\omega}\right)$ stands for the closure of ${\mathcal D}(\Omega)$ in $W^{1,p}\left(\Omega,\frac{d^n}{\omega}\right)$.
\begin{theo} \label{sobolev}
Let $p\in (1,+\infty)$. Then, for all $f\in L^p_0(\Omega)$, there exists $u\in W^{1,p}_0(\Omega,\frac{d^n}{\omega})$ such that
$$
\mbox{div }u=f\in \Omega
$$
and
\begin{equation} \label{sobestim}
\left\Vert Du\right\Vert_{L^p\left(\Omega,\frac{d^n}{\omega}\right)}\leq C\left\Vert f\right\Vert_{L^p(\Omega)},
\end{equation}
where the constant $C>0$ only depends on $\Omega$ and $p$. 
\end{theo}
In the statement of Theorem \ref{sobolev}, by $Du$, we mean the matrix $\left(\partial_iu_j\right)_{1\leq i,j\leq n}$, and the estimate (\ref{sobestim}) means that 
$$
\sum_{1\leq i,j\leq n} \left\Vert \partial_iu_j\right\Vert_{L^p\left(\Omega,\frac{d^n}{\omega}\right)}\leq C\left\Vert f\right\Vert_{L^p(\Omega)}.
$$
The proof of Theorem \ref{sobestim} goes through a decomposition of functions in $L^p_0(\Omega)$, interesting in itself, which follows from the solvability of the divergence problem established in the previous sections:
\begin{pro} \label{atomic}
Consider the Whitney decomposition of $\Omega$ already used in Section \ref{ineqlinfty}. Let $p\in (1,+\infty)$ and $f\in L^p_0(\Omega)$. Then, the function $f$ can be decomposed as
\begin{equation} \label{decompo}
f=\sum_{j} f_j,
\end{equation}
where, for all $j$, the function $f_j$ satisfies the following three properties:
\begin{itemize}
\item[$(i)$]
$f_j$ is supported in $2Q_j$,
\item[$(ii)$]
$\int_{2Q_j} f_j(x)dx=0$,
\item[$(iii)$]
$\sum_j \int_{2Q_j} \left\vert f_j(x)\right\vert^p\left(\frac{d^n(x)}{\omega(x)}\right)^pdx\leq C\int_{\Omega} \left\vert f(x)\right\vert^pdx$,
\end{itemize} 
for some $C>0$ only depending on $\Omega$, $p$ and the choice of the paths $\gamma$. \par
\noindent If one furthermore assumes that there exists $C>0$ such that the weight $\omega$ satisfies:
\begin{equation} \label{hypomega}
\omega(x)\leq Cd^{n}(x)\mbox{ for almost all }x\in \Omega,
\end{equation}
then the conclusion of $(iii)$ can be strengthened as
\begin{equation} \label{fjf}
C^{-1}\int_{\Omega} \left\vert f(x)\right\vert^pdx\leq \sum_j \int_{2Q_j} \left\vert f_j(x)\right\vert^p dx\leq C\int_{\Omega} \left\vert f(x)\right\vert^pdx.
\end{equation}
\end{pro}
{\bf Proof: } Let $f\in  L^p(\Omega)$ with zero integral and, applying Theorem \ref{divp}, pick up a vector field $u$ such that $f=\mbox{div }u$ and $\left\Vert \frac{d^{n-1}}{\omega}u\right\Vert_{L^{p}(\Omega)}\leq C\left\Vert f\right\Vert_{L^p(\Omega)}$. Let $(\chi_j)_{j\geq 1}$ be a partition of unity associated with the $Q_j$'s, {\it i.e.} a sequence of ${\mathcal D}(\R^n)$ functions satisfying $\sum_j \chi_j=1$ on $\R^n$, such that, for all $j\geq 1$, $\chi_j=1$ on $Q_j$, is supported in $2Q_j$ and satisfies $\left\Vert \nabla \chi_j\right\Vert_{\infty} \leq Cl_j^{-1}$. Setting
\[
f_j:=\mbox{div }(\chi_ju)
\]
for all $j\geq 1$, one clearly has (\ref{decompo}). It follows from the support property of $\chi_j$ that $f_j$ is supported in $2Q_j$. That $f_j$ has zero integral on $\Omega$ follows by an integration by parts on $2Q_j$ and the fact that $\chi_j=0$ on $\partial Q_j$. Finally, since $\omega(x)\geq cd^n(x)$ for all $x\in \Omega$,
\[
\begin{array}{lll}
\displaystyle \sum_{j\geq 1} \int_{\Omega} \left\vert f_j(x)\right\vert^p\left(\frac{d^n(x)}{\omega(x)}\right)^pdx & \leq & \displaystyle C\sum_{j\geq 1} \int_{2Q_j} \left\vert \chi_j(x)\right\vert^p \left\vert f(x)\right\vert^p\left(\frac{d^n(x)}{\omega(x)}\right)^pdx \\
& + & \displaystyle  C\sum_{j\geq 1} \int_{2Q_j} \left\vert u(x)\right\vert^p \left\vert \nabla \chi_j(x)\right\vert^p \left(\frac{d^n(x)}{\omega(x)}\right)^pdx\\
& \leq & \displaystyle C\int_{\Omega} \left\vert f(x)\right\vert^pdx+ C\sum_{j\geq 1} l_j^{-p} \int_{2Q_j} \left\vert u(x)\right\vert^p\left(\frac{d^{n-1}(x)}{\omega(x)}\right)^pdx\\
& \leq & \displaystyle C\int_{\Omega} \left\vert f(x)\right\vert^pdx,
\end{array}
\]
which is $(iii)$.\par
\noindent Assume now that (\ref{hypomega}) holds. Then,
$$
c\leq \frac{d^n(x)}{\omega(x)}\leq C,
$$
and
$$
\int_{\Omega} \left\vert f(x)\right\vert^pdx  \leq C\sum_j \int_{2Q_j} \left\vert f_j(x)\right\vert^pdx,
$$
which completes the proof of Proposition \ref{atomic}. \hfill\fin\par
\begin{rem} An analogous decomposition for functions in $L^p(\Omega)$ with zero integral was established by Diening, Ruzicka and Schumacher in \cite{drs}. Actually, in their result, $f$ can be taken in $L^p(\Omega,w)$ where the weight $w$ belongs to the Muckenhoupt class $A_q$ for some $q\in (1,+\infty)$, and (\ref{fjf}) holds with the Lebesgue measure in both sides on the inequalities, but the authors assume that $\Omega$ satisfies an emanating chain condition (whereas no assumption is made on $\Omega$ in Proposition \ref{atomic}). The proof in \cite{drs} is direct, which means that it does not use the divergence operator, but, as in the present paper, some consequences are derived for the solvability of the divergence operator in Sobolev spaces.
\end{rem}
{\bf Proof of Theorem \ref{sobolev}: } let $f\in L^p(\Omega)$ with zero integral and decompose $f=\sum f_j$ as in Proposition \ref{atomic}. For each $j$, there exists $u_j\in W^{1,p}_0(2Q_j)$ such that $\mbox{div }u_j=f_j$ and $\left\Vert Du_j\right\Vert_{L^p(2Q_j)}\leq C\left\Vert f_j\right\Vert_{L^p(2Q_j)}$ (see for instance \cite{bb}, Section 7.1, Theorem 2). Define $u:=\sum_j u_j$.  One clearly has $\mbox{div }u=f$. Moreover, by construction, there exist $0<c<C$ with the following property: for $j\geq 1$, there exists $\omega_j\in (0,+\infty)$ such that, for all $x\in Q_j$,
$$
c\omega_j\leq \omega(x)\leq C\omega_j.
$$
As a consequence, 
$$
\begin{array}{lll}
\displaystyle \int_{\Omega} \left\vert Du(x)\right\vert^p \left(\frac{d^n(x)}{\omega(x)}\right)^pdx & \leq & \displaystyle C\sum_j \int_{2Q_j} \left\vert Du_j(x)\right\vert^p \left(\frac{d^n(x)}{\omega(x)}\right)^pdx\\
& \leq  & \displaystyle  C\sum_j \frac{l_j^{np}}{\omega_j^p} \int_{2Q_j} \left\vert Du_j(x)\right\vert^p dx\\
& \leq & \displaystyle C\sum_j \frac{l_j^{np}}{\omega_j^p} \int_{2Q_j} \left\vert f_j(x)\right\vert^p dx\\
& \leq & \displaystyle C\int_{2Q_j} \left\vert f_j(x)\right\vert^p \left(\frac{d^n(x)}{\omega(x)}\right)^pdx\\
& \leq & \displaystyle C\int_{\Omega} \left\vert f(x)\right\vert^pdx. \hfill\fin
\end{array}
$$
As a corollary of Theorem \ref{sobolev}, we obtain:
\begin{cor}
Let $p\in (1,n)$ and set $p^{\ast}:=\frac{np}{n-p}$. Then, for all $f\in L^p_0(\Omega)$, there exists $u\in L^{p^{\ast}}\left(\Omega,\frac{d^n}{\omega}\right)$ such that
$$
\mbox{div }u=f\in \Omega
$$
and
\begin{equation} \label{sobestimbis}
\left\Vert u\right\Vert_{L^{p^{\ast}}\left(\Omega,\frac{d^n}{\omega}\right)}\leq C\left\Vert f\right\Vert_{L^p(\Omega)},
\end{equation}
where the constant $C>0$ only depends on $\Omega$ and $p$. 
\end{cor}
{\bf Proof: } consider again $f\in L^p(\Omega)$ with zero integral and let $f_j,u_j$ and $u$ be as in the proof of Theorem \ref{sobolev}. For all $j\geq 1$, by the usual Sobolev embedding,
$$
\int_{2Q_j} \left\vert u_j(x)\right\vert^{p^{\ast}}dx\leq C\left(\int_{2Q_j} \left\vert Du_j(x)\right\vert^pdx\right)^{\frac {p^{\ast}}p}\leq C\left(\int_{2Q_j} \left\vert f_j(x)\right\vert^pdx\right)^{\frac {p^{\ast}}p}.
$$
It follows that
$$
\begin{array}{lll}
\displaystyle \int_{\Omega}  \left\vert u(x)\right\vert^{p^{\ast}}\left(\frac{d^n(x)}{\omega(x)}\right)^
{p^{\ast}}dx & \leq & \displaystyle C\sum_{j} \int_{2Q_j}  \left\vert u_j(x)\right\vert^{p^{\ast}}\left(\frac{d^n(x)}{\omega(x)}\right)^{p^{\ast}}dx\\
& \leq & \displaystyle C\sum_j \frac{l_j^{np^{\ast}}}{\omega_j^{p^{\ast}}} \int_{2Q_j} \left\vert u_j(x)\right\vert^{p^{\ast}}dx\\
& \leq & \displaystyle C\sum_{j} \frac{l_j^{np^{\ast}}}{\omega_j^{p^{\ast}}} \left(\int_{2Q_j} \left\vert f_j(x)\right\vert^pdx\right)^{\frac {p^{\ast}}p}\\
& \leq & \displaystyle C\sum_j \left(\int_{2Q_j} \left\vert f_j(x)\right\vert^p\left(\frac{d^n(x)}{\omega(x)}\right)^{p^{\ast}}dx\right)^{\frac {p^{\ast}}p}\\
& \leq & \displaystyle C\left(\sum_{j} \int_{2Q_j} \left\vert f_j(x)\right\vert^p\left(\frac{d^n(x)}{\omega(x)}\right)^{p^{\ast}}dx\right)^{\frac {p^{\ast}}p}\\
& \leq & \displaystyle C\left(\int_{\Omega} \left\vert f(x)\right\vert^pdx\right)^{\frac{p^{\ast}}p},
\end{array}
$$
which ends the proof. \hfill\fin

\SE{Examples} \label{ex}

In this section we give various examples of domains $\Omega$ to which the preceding results apply and for which we obtain an estimate of the weight $\mbox{w}$. \par
\noindent In their paper \cite{adm}, the authors show how to invert the divergence for John domains, with the weight $\mbox{w}(x) = d(x)$ (note that the corresponding weighted Poincar\'e inequalities on John domains were established in \cite{dd}). Our method gives a generalization to the so-called $s$-John domains, that we will present first; then we will turn to a particular case which exhibits special features, that of strongly h\"olderian domains.

\subsection{The case of $s$-John domains}

Let $s \geq 1$ be a fixed parameter. Recall that $\Omega$ is an $s$-John domain when, $x_0$ being a chosen reference point, there exist a constant $C > 0$ and, for each $y \in \Omega$, a rectifiable path $\gamma(y)$ included in $\Omega$ which links $y$ to $x_0$ and satisfies the estimate
\begin{equation}\label{s}
\forall \tau \in [0, l(\gamma(y))] \  d(\gamma(\tau,y)) \geq C \tau^s,
\end{equation}
where $\tau$ is the arc-length and $l(\gamma(y))$ the total length of the path $\gamma(y)$. When $s=1$, $\Omega$ is a John domain. Using the same algorithm as in the proof of Theorem \ref{solvability}, we modify the initial family $\gamma$ so as to ensure properties $(\gamma.a)$ to $(\gamma.d)$. We let the reader check that property (\ref{s}) remains preserved. Here are explicit examples. \par

\medskip

\noindent {\bf Example 1} Let $\Omega$ be the complement in $B(0,10)$ of the logarithmic spiral 
$$
\Gamma = \{x=r(\cos\theta, \sin\theta); r=0 \mbox{ or } r \in ]0,1] \mbox{ with } r=e^{-\theta}\}.
$$
Then $\Omega$ is a John domain. \par
\medskip
\noindent {\bf Example 2} Let $\Omega$ be the complement in $B(0,10)$ of the spiral
$$
\Gamma = \{x=r(\cos\theta, \sin\theta); r=0 \mbox{ or } r \in ]0,1] \mbox{ with } \theta = r^{-a}\}
$$
for some $a>0$. Then $\Omega$ is an $s$-John domain if and only if $a<1$, with $s=\frac{1+a}{1-a}$. Note, however, that $\Omega$ fulfills (\ref{dOmega}) if and only if $a<3$. \par

\medskip

\noindent Returning to the general $s$-John domains, we now assume that $\Omega$ fulfills hypothesis (\ref{dOmega}). Applying Lemma \ref{solutioninfty} and the preceding sections, we can invert the divergence on various $L^p$ spaces on $\Omega$ with the weight
$$
\mbox{w}(x)=\omega(x) d(x)^{-n+1},
$$
where $\omega$ is defined by (\ref{defomega}). We then have:

\begin{lem}\label{estimw-sJohn}
$\omega(x) \leq C d(x)^{\frac{n}{s}}.$
\end{lem}

Hence we recover the result of \cite{hkjones}, Theorem 7, that on $s$-John domains the following Poincar\'e inequality holds:
$$
\int_{\Omega} \left\vert g(x)-g_{\Omega}\right\vert dx \leq C(\Omega) \int_{\Omega}  d(x)^{\frac{n}{s}-n+1}
\left\vert \nabla g(x)\right\vert dx.
$$
It is known, also, that $\frac{n}{s}-n+1$ is the best possible exponent in the class of $s$-John domains (see \cite{bstr}). Let us, however, remark that the estimate in this Lemma may be of poor quality, since it may happen for large values of $s$ that
$$
\int_0^{+\infty} d(x)^{\frac{n}{s}-n+1} dx = + \infty,
$$
while we know that $\mbox{w} \in L^{1}(\Omega)$. This reveals that, at least in some cases, the family of weights $d(x)^\alpha$, $\alpha \in \R$, is not rich enough to accurately describe what happens. \par

\medskip

\noindent The proof of Lemma \ref{estimw-sJohn} is elementary. Let $x,y \in \Omega$ such that $\mbox{dist}_{\Omega}(\gamma(y),x) \leq \frac{1}{2} d(x)$. There exists $\tau \in [0, l(\gamma(y))]$ verifying
$$
\left\vert \gamma(\tau,y) -x \right\vert \leq \frac{1}{2} d(x),
$$
and thus
$$
\tau^s \leq C d(\gamma(\tau,y)) \leq C d(x).
$$
Since $\left\vert y-x \right\vert \leq \tau + \frac{1}{2} d(x)$, we have $y \in B(x,C d(x)^{\frac1{s}})$, which implies
$\omega(x) \leq C d(x)^{\frac{n}{s}}$ as desired.

\subsection{The case of strongly h\"olderian domains}

There is, however, a particular case where the above estimate can be improved as soon as $s>1$. This is when $\Omega$, in addition to the preceding hypothesis, fulfills the requirements of the following definition.

\begin{defi}
Say that $\Omega$ is a strongly h\"olderian domain when there exist $\alpha \in ]0,1]$ and an integer $N\geq 1$, $N$ functions $\Phi_j \in C^{\alpha}(\R^{n-1})$, $N+1$ domains $O_j$ in $\R^n$ and $N$ isometries $\eta_j$ of $\R^n$ such that
    \begin{itemize}
    \item[\bf{$\cdot$}] $\Omega \subset \bigcup_{j=0}^{N}O_j$,
    \item[\bf{$\cdot$}] $\overline{O_0} \subset \Omega$,
    \item[\bf{$\cdot$}] when $x \in O_j$, $j \geq 1$, then 
    $x \in \Omega \mbox{ if and only if } \tilde{x_n} < \Phi_j(\tilde{x'}),$
    where $\eta_j (x) = (\tilde{x'},\tilde{x_n}), \tilde{x'} \in \R^{n-1}, \tilde{x_n} \in \R$,
    \item[\bf{$\cdot$}] there exists $h>0$ such that, if $x\in O_j \cap \Omega, j \geq 1$, and $\tilde{x_n} < \Phi_j (\tilde{x'})-h$, then $x \in O_0$,
    \item[\bf{$\cdot$}] conversely, for any $x\in O_j \cap \Omega, j \geq 1$, there exists $z \in O_0$ such that $\tilde{z'} = \tilde{x'}$ and $\tilde{z_n} \leq \tilde{x_n}$, where $\eta_j (z) = (\tilde{z'},\tilde{z_n}).$
    \end{itemize}
\end{defi}

\noindent We construct a family $\gamma$ of paths satisfying properties ($\gamma.a$) to ($\gamma.d$) and adapted to the nature of the boundary of $\Omega$.\par
\noindent Choose $x_0 \in O_0$ and let $y \in \Omega$. When $y \in O_0$, we cover $O_0$ by a finite number of Whitney cubes associated to $\Omega$ and apply the algorithm described in the proof of Theorem \ref{solvability} to obtain a family of paths $\gamma(y)$ which lie inside $O_0$.\par
\noindent When $y \in O_j, j \geq 1$, we assume for simplicity that $\eta_j = Id$ (which is possible by rotating $\Omega$). Write $y=(y',y_n)$ and select $z=(y',z_n) \in O_0$ with $z_n \leq y_n$. We link $y$ to $z$ by the vertical segment $\{ (y',(1-t)y_n + tz_n); t \in [0,1]\}$, and then $z$ to $x_0$ as described above. \par
\noindent The reader may check that the resulting family of paths satisfies properties ($\gamma.a$) to ($\gamma.d$), and $\Omega$ is an $s$-John domain for $s=\frac{1}{\alpha}$.

\medskip

\noindent Let $\mbox{w}$ be the associated weight, which allows to solve the problem $(\div)_{\infty}$ in $\Omega$ thanks to Lemma \ref{solutioninfty}. We have the following improvement over the result of Lemma \ref{estimw-sJohn}:

\begin{lem}
There exists a constant $C>0$ such that  
$$ \forall x \in \Omega, \, \, \mbox{w}(x) \leq C d(x)^\alpha.$$
\end{lem}

\noindent Since $d^\alpha$ is always integrable over $\Omega$, this Lemma implies that, for this class of domains, the problem $(\div)_{\infty}$ is solvable with the weight $d^\alpha.$ \par

\medskip

\noindent Let us turn to the proof and take $x,y \in \Omega$ with $\mbox{dist}_{\Omega}(\gamma(y),x) \leq \frac12 d(x)$. We may assume that $d(x) < \delta$ for some $\delta > 0$ to be chosen, since otherwise there is nothing to prove. If $\delta$ is small enough, there exists $j \geq 1$ such that $x,y \in O_j$. We again simplify the argument by assuming $\eta_j = Id.$ Write $y=(y',y_n)$ and $z=(y',z_n)$ as above. The distance from $\gamma(y)$ to $x$ is not attained at a point lying in $\overline O_0$ provided $\delta$ is sufficiently small, so that there exists $t \in [0,1[$ satisfying
$$
\left\vert (y', (1-t)y_n + tz_n) - x \right\vert \leq \frac12 d(x).
$$
Hence we have
\begin{equation}\label{star}
\left\vert y' - x' \right\vert \leq \frac12 d(x)
\end{equation}
and (recall that $z_n < y_n$)
\begin{equation}\label{b}
x_n - \frac12 d(x) \leq y_n < \Phi_j (y').
\end{equation}
We want to show that
\begin{equation}\label{dstar}
\left\vert y_n - x_n \right\vert \leq C d(x)^\alpha. 
\end{equation}
From (\ref{star}) and the $\alpha$-regularity of $\Phi_j$, we deduce
$$
\left\vert \Phi_j (y') - \Phi_j (x') \right\vert \leq C d(x)^\alpha
$$
and
\begin{equation}\label{c}
y_n < \Phi_j (x') + C d(x)^\alpha.
\end{equation}
Let $(u',\Phi_j(u'))$ be a point of $\partial \Omega$ at which $d(x)$ is attained (such a point exists once $\delta$ is small enough and $j$ is appropriately chosen). We have $\left\vert x' - u' \right\vert \leq d(x)$ and $\left\vert \Phi_j(u') - x_n \right\vert \leq d(x),$ whence 
$$
\begin{array}{lll}
\left\vert \Phi_j(x') - x_n \right\vert & \leq & \left\vert \Phi_j(x') - \Phi_j(u') \right\vert + d(x) \\
                                        & \leq & C d(x)^\alpha.
\end{array}
$$
From (\ref{b}) and (\ref{c}) we then have
$$
x_n - \frac12 d(x) \leq y_n \leq x_n + C d(x)^\alpha,
$$
which proves (\ref{dstar}). With (\ref{star}), this implies
$$
\omega (x) \leq C d(x)^{n-1+\alpha},
$$
and the proof is readily ended.

\end{document}